\documentclass[11pt]{article}
\usepackage{amsmath, amsfonts, amsthm, amssymb, amscd}
\usepackage[mathcal]{euscript}
\usepackage{dsfont, pxfonts, mathrsfs, accents, verbatim}
\theoremstyle{plain}
        \newtheorem{theorem}{Theorem}[section]

        \newtheorem{proposition}[theorem]{Proposition}
        \newtheorem{lemma}[theorem]{Lemma}
        \newtheorem{corollary}[theorem]{Corollary}

\numberwithin{equation}{section}

\newcommand \be     {\begin{equation}}
\newcommand \ee     {\end{equation}}

\newcommand \del        \partial
\newcommand \eps        \varepsilon
\newcommand \auth   \textsc

\newenvironment{pf}{{\noindent \it  Proof.}}{{\hfill q.e.d.}\\}

\begin{document}

\title{Strong Uniqueness of the Ricci Flow}
\author{Bing-Long Chen\footnote{
Department of Mathematics, Sun Yat-Sen University 510275, Guang
Zhou, P.R.of China. E-mail: {\sl mcscbl@mail.sysu.edu.cn.}
\newline
}}
\maketitle

\begin{abstract}
In this paper,  we  derive some local a priori estimates for Ricci
flow. This gives rise to some  strong uniqueness theorems.  As a
corollary, let $g(t)$ be a smooth complete solution to the Ricci
flow on $\mathbb{R}^{3}$,  with the canonical Euclidean metric $E$
as initial data,  then $g(t)$ is trivial, i.e. $g(t)\equiv E$.
\end{abstract}


 \section{Introduction}
\vskip 0.5cm

\qquad  The Ricci flow $ \frac{\partial}{\partial
t}g_{ij}(x,t)=-2R_{ij}(x,t)$,  was introduced by Hamilton in
\cite{Ha1}. The major application of this equation to  lower
dimensional topology has had a great impact in modern mathematics
(see \cite{Ha1}, \cite{Ha2}, \cite{Ha3}  \cite{P1}, \cite{P2}).
The power of these geometric applications grew out of the
fundamental PDE theory of the equation.  These two aspects had
been intertwined all the time since the foundation of the Ricci
flow.

In this paper, we  go back to  some fundamental  PDE problems of
this equation.

 Let's look at one heuristic analogue, the standard
heat equation $(\frac{\partial}{\partial t}-\triangle)u=0$ on
$\mathbb{R}^{n}.$  If $u$ grows  slower than function
$e^{a|x|^{2}}$ for some $a>0,$ then $u$ is unique for all such
solutions with same initial data. Moreover, if $|u|\mid_{t=0}\leq
Ce^{a|x|^{2}},$ it is not hard to see the short time existence (of
solutions of same type) from the heat kernel convolution.   For
Ricci flow, the Ricci curvature behaves like twice derivative of
logarithmic of the metric. So bounded curvature condition for
Ricci flow resembles  growth $e^{a|x|^{2}}$ for standard heat
equation. Actually, the fundamental work \cite{Sh1} showed that on
complete manifolds with bounded curvature, the Ricci flow always
admits  short time
 solutions of bounded curvature.  X.P.Zhu and the author  recently
\cite{CZ} proved that the uniqueness theorem holds for solutions in
the class of bounded curvature. For an interesting application of
this theorem to the theory of Ricci flow with surgery, we refer the
readers to see \cite{CZ05F} or  relevant discussions in
\cite{CaoZ}\cite{KL}\cite{MT}.

 However,  if one don't impose
any growth conditions, the solutions to the heat equation
$(\frac{\partial}{\partial t}-\triangle)u=0$  are no longer unique.
 For instance, when $n=1,$ the famous Tychonoff's example
$u(x,t)=\sum_{k=0}^{\infty}\frac{x^{2k}}{(2k)!}\frac{d^{k}}{dt^{k}}e^{-\frac{1}{t^{2}}},$
is a smooth nontrivial solution to the heat equation with $0$
initial data. The purpose of this paper is to investigate the
analogous problem for Ricci flow. Nevertheless, Ricci flow,
 as the most natural intrinsic heat deformation of metrics, has quite complicated nonlinearity.
 We attempt to show that, in
certain extent, the above phenomenon never
 happens for geometrically  reasonable solutions.

Now, we formulate one of the main results of this paper as
following

\begin{theorem}
\label{T1.1} Let $(M,g(x))$ be a complete noncompact three
dimensional manifold with bounded and nonnegative sectional
curvature $
 0\leq Rm\leq K_0,
$ for some fixed constant $K_0.$ Let ${g}_{1}(x,t),$
${g}_{2}(x,t),$ $t\in[0,T],$ be two smooth complete solutions to
the Ricci flow with initial data $g(x).$ Then we have
$g_{1}(t)\equiv g_2(t),$ for $0\leq t < \min\{T,\frac{1}{4K_0}\}.$
\end{theorem}
A simple example is the Euclidean space $ \mathbb{R}^{3}:$
\begin{corollary}
\label{C1.2} Let $g(t),$  $t\in [0,T]$, be a smooth complete
solution to the Ricci flow on $ \mathbb{R}^{3},$  starting with
the canonical Euclidean metric $E$, then $g(t)\equiv E.$
\end{corollary}
The most important feature of these uniqueness theorems,
   is that we do not require any extra growth conditions on the
   solutions
  except the geodesic completeness.

Theorem \ref{T1.1} may be viewed as a generalization of \cite{CZ},
and we call it a strong uniqueness theorem (see the extrinsic
version in \cite{CY}).

  In views of \cite{CZ}, the whole issue is reduced to the curvature estimates.
  In this paper,  we will derive some local
curvature estimates in dimension 3 (or 2),  which have  their own
interest from PDE point of view. Our strategy is the following.
 The singularities of ancient type
occur naturally, once the desired estimate fails. Through the
great works of Hamilton and Perelman, the structures of
singularities in dimension 3 have already been well-understood
nowadays. One crucial reason why Ricci flow works in dimension $3$
in the classical theory is that we have Hamilton-Ivey's curvature
pinching estimate, which guarantees the singularities are always
nonnegative curved. Recall these estimates were proved by maximum
principle.  In the classical setting, this principle can only be
applied on manifolds or (Ricci flow) solutions with bounded or
suitable  growth curvature. Remember the curvature estimate is
just the goal we want to achieve. To go around this difficulty, in
this paper, we will derive some pinching estimates of similar
type, but in a purely local way (see section 2)

In this regard, let us recall  the so-called pseudolocality
theorem of Perelman \cite{P1}. The point we should mention here,
is that the pseudolocality theorem of Perelman \cite{P1} is
basically proved for compact manifolds.  Since the justification
of integration by parts on the whole manifold is also ultimately
related to the geometry of the solution, this makes the situation
very complicated(see the proof in section 10 in \cite{P1}).
 Actually, it is a still an open problem if
the pseudolocality theorem holds for \textbf{any} complete
solutions to the Ricci flow. Recently, in \cite{CTY}, by assuming
the solution has bounded curvature, the pseudolocality theorem of
Perelman has been generalized to complete manifolds. As mentioned
above, the key point for our strong uniqueness is just the
curvature bound.  In this paper, we will adopt a totally different
approach.

We remark that in dimension 2, we even have a  better strong
uniqueness theorem, i.e. nonnegative curvature assumption can be
removed (see Theorem \ref{T3.10}).

 The paper
is organized as follows. In section $2$, we derive a pure local
pinching estimate for $3$ dimensional Ricci flow. In section $3$,
we will show various local a priori curvature estimates, which may
give rise to the proof of the uniqueness theorems. In section $4$,
we will discuss some further open problems.

\textbf{Acknowledgements} The author is grateful to  Professor
X.P.Zhu for many stimulating discussions. The author is also
thankful to Professors D. Burago,  X. C.  Rong and some others for
the discussions on the injectivity radius assumption in the
previous version of the paper.  Finally, M. Simon told the author
that he had similar interior estimates (provided the curvature of
the solution has been controlled by $\frac{K}{t}$) as in our
Theorem \ref{T3.1} in his habilitation thesis \cite{Si}.  This
work was partially supported by grants from Sun Yat-Sen University
(NCET-050717,  N0. 34000-3171404 $\&$ 34000-1131040).


\section{Local pinching estimate}
 Hamilton-Ivey's pinching estimate plays a
substantial role in the application of the Ricci flow to the
geometrization conjecture in dimension $3$. As a matter of fact,
this is one main reason why this theory works in this very
dimension. As we mentioned in the introduction, the proof of this
estimate is by maximum principle for compact or complete solutions
with bounded curvature. Because the curvature bound on the whole
manifold is just the goal we want to achieve, this becomes an
obstacle for us. Fortunately, we find that the equation has
certain good nonlinearity, which enables us to localize all the
estimates.

We start with the local estimate of scalar curvature, which is
dimensionally free.
\begin{proposition}
\label{P2.1} For any $0<\delta<\frac{2}{n}$ , there is
$C=C(\delta,n)>0$ satisfying the following property. Suppose we have
a smooth solution $g_{ij}(x,t)$ to the Ricci flow on an $n$
dimensional manifold $M$, such that for any $t\in [0,T]$,
 $B_{t}(x_0,Ar_0)$ are compactly
contained in $M$ and assume that
 $Ric(x,t)\leq ({n-1}){r_0}^{-2}$
for $x\in B_{t}(x_0,r_0),$ $t\in [0,T]$ and   $R\geq -K$ ($K\geq
0$) on $B_{0}(x_0,Ar_0)$ at $t=0.$ Then we have

(i) $ R(x,t)\geq \min\{-\frac{1}{(\frac{2}{n}-\delta)t+\frac{1}{K}},
-\frac{C}{Ar_0^{2}}\},$ if $A\geq2;$

(ii) $ R(x,t)\geq
\min\{-\frac{1}{(\frac{2}{n}-\delta)t+\frac{n}{K}},
-\frac{C}{A^{2}r_0^{2}}\}, $ if  $A\geq
\frac{40}{3}(n-1)r_0^{-2}T+2,$

 whenever $x\in B_{t}(x_0,\frac{3A}{4}r_0)$, $t\in[0,T].$
\end{proposition}
 \begin{pf}  By \cite{P1}, we have
\begin{equation}
\label{2.1} (\frac{\partial}{\partial
t}-\triangle)d_{t}(x_0,x)\geq -\frac{5(n-1)}{3}r_0^{-1},
\end{equation} whenever
$d_{t}(x,x_0)>r_0,$ in the sense of support functions.

We divide the discussion into two cases.

Case(a): $A\geq \frac{40}{3}(n-1)r_0^{-2}T+2.$

We consider the function
$u=\varphi(\frac{d_{t}(x_0,\cdot)+\frac{5(n-1)r_0^{-1}t}{3}}{Ar_0})R,$
where $\varphi$ is a fixed smooth nonnegative non-increasing
function such that $\varphi=1$ on $(-\infty, \frac{7}{8}],$  and
$\varphi=0$ on $[1,\infty).$

It is clear
\begin{equation}
\label{2.2}
\begin{split}
(\frac{\partial}{\partial t}-\triangle)u=& \varphi'R
\frac{1}{Ar_0}[(\frac{\partial}{\partial
t}-\triangle)d_{t}(x_0,x)+\frac{5}{3}(n-1)r_0^{-1}]\\
 & -\varphi''\frac{1}{(Ar_0)^{2}}R+2\varphi |Ric|^{2}-2\nabla
\varphi\cdot\nabla R,
\end{split}
\end{equation}
at smooth points of distance function.

 Let $u_{min}(t)=\min_{M}
u(\cdot,t)$. If $u_{min}(t_0)\leq0$ and $u_{min}(t_0)$ is achieved
at some point $x_1,$ then $\varphi 'R(x_1,t_0)\geq 0.$ Hence, by
(\ref{2.1}), the first term in the right hand side of (\ref{2.2}) is
nonnegative.  Now by applying the maximum principle and standard
support function technique,
 we have for any small $\delta>0$
\begin{equation}
\begin{split}
\label{e2.1.1} \frac{d^{-}}{d
t}u_{min}\mid_{t=t_0}&:=\liminf_{\triangle t\searrow
0}\frac{u_{min}(t_0+\triangle t)-u_{min}(t_0)}{\triangle t}\\
&\geq \frac{2}{n}\varphi
R^{2}+\frac{1}{(Ar_0)^{2}}(\frac{2\varphi'^{2}}{\varphi}-\varphi'')R\\&
\geq(\frac{2}{n}
-\delta)u_{min}(t_0)^{2}+\frac{\delta}{2}(u_{min}(t_0)^{2}-\frac{C^{2}}{(Ar_0)^{4}}).
\end{split}
\end{equation}
 provided $u_{min}(t_0)\leq0,$  where we have used $|
 \frac{2\varphi'^{2}}{\varphi}-\varphi''|\leq C\sqrt{\varphi}$ and
 Cauchy-Schwartz
 $|\frac{1}{(Ar_0)^{2}}(\frac{2\varphi'^{2}}{\varphi}-\varphi'')R|
 \leq \frac{\delta}{2} \varphi R^{2}+\frac{C}{(Ar_0)^{4}}.$

 By integrating the  inequality
(\ref{e2.1.1}), we get
$$
u_{min}(t)\geq \min\{-\frac{1}{(\frac{2}{n}-\delta)t+\frac{1}{K}},
-\frac{C}{(Ar_0)^{2}}\}.
$$
This implies
$$
R(x,t)\geq \min\{-\frac{1}{(\frac{2}{n}-\delta)t+\frac{1}{K}},
-\frac{C(\delta)}{(Ar_0)^{2}}\},
$$
 whenever $x\in B_{t}(x_0, \frac{3A}{4}r_0).$

Case(b):  $A\leq \frac{40}{3}(n-1)Tr_0^{-2}+2.$

Consider the function $u=\varphi(\frac{d_{t}(x_0,\cdot)}{Ar_0})R,$
the similar argument yields
$$
 u_{min}(t)\geq \min\{-\frac{1}{(\frac{2}{n}-\delta)t+\frac{1}{K}},
-\frac{C(\delta)}{Ar_0^{2}}\}.
$$
The proof is completed.
\end{pf}

In dimension 3,  in terms of  moving frames \cite{Ha2}, the
curvature operator, $M_{ij}=Rg_{ij}-2R_{ij}$,   has the following
evolution equation
$$
\frac{\partial}{\partial t}M=\triangle M+M^{2}+M^{\sharp},
$$
 where $M^{\sharp}$ is
the lie algebra adjoint of $M.$ Let $\lambda\geq \mu\geq \nu$ be
the eigenvalues of $M$, the same eigenvectors also diagnolize
$M^{2}+M^{\sharp}$ with eigenvalues $\lambda^{2}+\mu\nu,$
$\mu+\lambda\nu,$ $\nu+\lambda\mu.$ The following estimate may be
viewed as a local version of Hamilton-Ivey pinching estimate.

\begin{proposition}
\label{P2.2} For any $k\in \mathbb{Z}_{+},$ there is $C_{k}$
depending only on $k$ satisfying the following property. Suppose we
have a smooth solution $g_{ij}(x,t)$ to the Ricci flow on a three
manifold $M$, such that for any $t\in [0,T]$, $B_{t}(x_0,Ar_0)$ are
compactly contained in $M$ and assume that
  $Ric(x,t)\leq (n-1)r_0^{-2}$for $x\in B_{t}(x_0,r_0),$ $t\in [0,T];$ and
 $\lambda+\mu+k\nu\geq
 -K_k$($K_k\geq 0$)
 on $B_0(x_0,Ar_0)$ at time $0$. Then we have

(i)  $ \lambda+\mu+k\nu\geq \min\{-\frac{C_{k}}{t+\frac{1}{K_k}},
-\frac{C_{k}}{Ar_0^{2}}\}, $ if $A\geq2;$

(ii) $ \lambda+\mu+k\nu\geq \min\{-\frac{C_{k}}{t+\frac{1}{K_k}},
-\frac{C_{k}}{A^{2}r_0^{2}}\}, $ if $A\geq
 \frac{40(n-1)k}{3} r_0^{-2}T+2,$

 whenever $x\in B_{t}(x_0,\frac{A}{2}r_0)$, $t\in[0,T].$
\end{proposition}
\begin{pf} We only prove the general case (i).   We will argue by induction on $k\in
\mathbb{Z}_{+}$ to prove the estimate holds on ball of radius
$(\frac{1}{2}+\frac{1}{2^{k_0+1}})Ar_0.$ The $k=1$ case follows
from Proposition \ref{P2.1}, and radius of the ball is
$\frac{3A}{4}r_0.$  Suppose we have proved the result for
$k=k_0\in \mathbb{Z}_{+},$ that is to say, there is constant
$C_{k_0}$ such that \be \label{2.4} \lambda+\mu+k_0\nu\geq
\min\{-\frac{C_{k_0}}{t+\frac{1}{K_{k_0}}},
-\frac{C_{k_0}}{Ar_0^{2}}\}, \ee whenever $x\in
B_{t}(x_0,(\frac{1}{2}+\frac{1}{2^{k_0+1}})Ar_0)$, $t\in[0,T].$ We
are going to prove the result for $k=k_0+1$ on ball of radius
$(\frac{1}{2}+\frac{1}{2^{k_0+2}})Ar_0.$

Without loss of generality, we may assume $K_1\leq K_2\leq K_3\leq
\cdots.$

Define a function $
C_{k_0}(t):=\max\{\frac{C_{k_0}}{t+\frac{1}{K_{k_0}}},
\frac{C_{k_0}}{Ar_0^{2}}\}.$

 Let
$$N_{ij}=Rg_{ij}+k_0M_{ij},\ \ \ \ \ \  P_{ij}=\varphi(\frac{d_t(x,x_0)}{Ar_0})
(Rg_{ij}+k_0M_{ij}),$$ where $\varphi$ is a smooth nonnegative
decreasing function, which is $1$ on $(-\infty,
\frac{1}{2}+\frac{1}{2^{k_0+2}}]$ and $0$ on
$[\frac{1}{2}+\frac{1}{2^{k_0+1}},\infty)$.  Note that the least
eigenvalue of $N_{ij}$ is $\lambda+\mu+(k_0+1)\nu.$ Let $V$ be the
corresponding (time dependent) unit eigenvector  of $N_{ij}$.

By direct computation, we have

$$
(\frac{\partial}{\partial
t}-\triangle)P_{ij}=-2\nabla_l\varphi\nabla_l N_{ij}+Q_{ij}
$$
where $Q_{ij}$ satisfies

 \begin{equation*}
\begin{split}
Q(V,V)&=\varphi(\lambda^{2}+\mu^{2}+(k_0+1)\nu^{2}+\mu\nu+\lambda\nu+(k_0+1)\lambda\mu)\\
& \ \ +[\varphi' \frac{1}{Ar_0}[(\frac{\partial}{\partial
t}-\triangle) d_{t}(x_0,x)]
 -\varphi''\frac{1}{(Ar_0)^{2}}](\lambda+\mu+(k_0+1)\nu).
\end{split}
\end{equation*}

Let $$u(t):=\min_{x\in M}(\lambda+\mu+(k_0+1)\nu)\varphi(x,t).$$

For fixed $t_0\in [0,T],$  assume
$(\lambda+\mu+(k_0+1)\nu)\varphi(x_0',t_0)=u(t_0)<-2C_{k_0}(t_0).$
Otherwise, we have the estimate at time $t_0.$

Combining with (\ref{2.4}), we have
$(\lambda+\mu+(k_0-1)\nu)(x_0',t_0)\geq 0.$ Note that
$\nu(x_0',t_0)$ is negative, otherwise
$(\lambda+\mu+(k_0+1)\nu)(x_0',t_0)\geq 0.$ Hence
$(\lambda+\mu)\varphi(x_0',t_0)\geq 0.$

 We compute
\begin{equation*}
\begin{split}
Q(V,V)(x_0',t_0)&=\varphi(\lambda^{2}+\mu^{2}+(k_0+1)\nu^{2}+(\lambda+\mu)\nu+(k_0+1)\lambda\mu)\\
& \ \ +[\frac{\varphi'}{\varphi}
\frac{1}{Ar_0}(\frac{\partial}{\partial t}-\triangle) d_{t}(x_0,x)
 -\frac{\varphi''}{\varphi}\frac{1}{(Ar_0)^{2}}]u(t_0)\\
 &=\varphi\frac{(\lambda+\mu+(k_0+1)\nu)^{2}}{(k_0+1)}-\frac{\lambda+\mu}{k_0+1}\varphi(\lambda+\mu+(k_0+1)\nu)
 \\
 &  \ \ \ +\varphi(\lambda^{2}+\mu^{2}+(k_0+1)\lambda\mu) +[\frac{\varphi'}{\varphi}
\frac{1}{Ar_0}(\frac{\partial}{\partial t}-\triangle) d_{t}(x_0,x)
 -\frac{\varphi''}{\varphi}\frac{1}{(Ar_0)^{2}}]u(t_0)\\
 & = I+II+III+IV.
\end{split}
\end{equation*}

Since $(\lambda+\mu)\varphi(x_0',t_0)\geq 0$ and $u(t_0)<0,$ we
have $II\geq 0.$ To deal with term $III,$ we divide into two
cases.

Case $(\alpha)$:
$\mu(x_0',t_0)<-\frac{\lambda(x_0',t_0)}{k_0+1}.$

By (\ref{2.4}),  $(\lambda+\mu+k_0\nu)(x_0,t_0)\geq-C_{k_0}(t_0)$,
we have $-\nu(x_0',t_0)\leq
\frac{\lambda(x_0',t_0)}{k_0+1}+\frac{C_{k_0}(t_0)}{k_0}.$ Hence
at $(x_0',t_0)$, we have
\begin{equation*}
\begin{split} \lambda^{2}+\mu^{2}+(k_0+1)\lambda\mu& \geq
\lambda^{2}+(\frac{\lambda}{k_0+1})^{2}-(k_0+1)\lambda(\frac{\lambda}{k_0+1}+\frac{C_{k_0}(t_0)}{k_0})\\
& \geq
(\frac{\lambda}{k_0+1})^{2}-(k_0+1)\frac{C_{k_0}(t_0)}{k_0}\lambda\geq
-\frac{(k_0+1)^{4}C_{k_0}(t_0)^{2}}{4k_0^{2}}.
\end{split}\end{equation*}

Case $(\beta)$:
$\mu(x_0',t_0)\geq-\frac{\lambda(x_0',t_0)}{k_0+1}.$

In this case,
$(\lambda^{2}+\mu^{2}+(k_0+1)\lambda\mu)(x_0',t_0)\geq 0$ holds
trivially.

Hence in either case, we have
$$\lambda^{2}+\mu^{2}+(k_0+1)\lambda\mu\geq -\frac{(k_0+1)^{4}C_{k_0}^{2}(t_0)}{4k_0^{2}}.$$
Therefore,
\begin{equation*}
\begin{split}
Q(V,V)(x_0',t_0)
 &\geq
 \varphi\frac{(\lambda+\mu+(k_0+1)\nu)^{2}}{(k_0+1)}-\frac{(k_0+1)^{4}C_{k_0}(t_0)^{2}}{4k_0^{2}}\varphi\\
  &  +[\frac{\varphi'}{\varphi}
\frac{1}{Ar_0}[(\frac{\partial}{\partial t}-\triangle)
d_{t}(x_0,x)]
 -\frac{\varphi''}{\varphi}\frac{1}{A^{2}r_0^{2}}]u(t_0)\\
& \geq
 \frac{1}{(k_0+1)\varphi}[u^{2}-(\frac{5(n-1)\varphi'}{3Ar_0^{2}}+\frac{k_0+1}
 {A^{2}r_0^{2}}\varphi'')u]-\frac{(k_0+1)^{4}C_{k_0}^{2}(t_0)}{4k_0^{2}}.
\end{split}
\end{equation*}
Since $|\varphi'|\leq C 2^{k_0},
|\varphi''|+\frac{\varphi'^{2}}{\varphi}\leq C 2^{2k_0},$  by
applying maximum principle, we have
\begin{equation*}
\begin{split}
\frac{d^{-}}{dt}\mid_{t=t_0}u& \geq
Q(V,V)(x_0',t_0)+\frac{2}{(Ar_0)^{2}}\frac{\varphi'^{2}}{\varphi^{2}}u(t_0)\\
 &\geq \frac{1}{2(k_0+1)}u^{2}
\end{split}
\end{equation*}
provided $|u|(t_0)\geq \max\{ CC_{k_0}(t_0){k_0}^{\frac{3}{2}},
C\frac{2^{2k_0}k_0}{Ar_0^{2}}\},$
 where $C$ is some universal
constant. By integrating the above differential inequality, we get
estimate:
$$u(t)\geq
\min\{\frac{1}{\frac{1}{u(0)}-\frac{t}{2(k_0+1)}},-CC_{k_0}(t){k_0}^{\frac{3}{2}},
-C\frac{2^{2k_0}k_0}{Ar_0^{2}}\}.$$

By the definition of $C_{k_0}(t),$ noting $-K_{k_0}\geq -K_{k_0+1}$,
clearly, there is a $C_{k_0+1}$ such that $$u(t)\geq
\min\{-\frac{C_{k_0+1}}{t+\frac{1}{K_{k_0+1}}},
-\frac{C_{k_{0}+1}}{Ar_0^{2}}\}.$$

The proof of case (ii) is similar. We use cut-off function
$\varphi(\frac{d_{t}(x_0,\cdot)+\frac{5(n-1)r_0^{-1}t}{3}}{Ar_0}),$
where $\varphi$ is a suitably chosen function which depends on
$k_0$ in the inductive step.

\end{pf}
We remark that by following the constants in the proof, the constant
$C_k$ may be chosen to be $Ck^{Ck}$ for some universal constant $C.$
The factor $\frac{1}{2}$ in the radius $\frac{1}{2}Ar_0$ is not
important, it may be replaced by any constant in $(0,1).$
\begin{corollary}
\label{C2.3} Suppose we have a complete smooth solution
$g_{ij}(x,t)$ to the Ricci flow on $M\times[0,T]$, then whenever
$t\in [0,T]$ we have

(i) if $R\geq -K$ for $0\leq K\leq \infty$ at $t=0,$ then
$$
R(\cdot,t)\geq -\frac{n}{2t+\frac{n}{K}};
$$

(ii) if dim $M$=3, then for any $k>0$, there is $C_k>0$ depending
only on $k$ such that if at $t=0$, $\lambda+\mu+k\nu\geq -K_k$ for
some $0\leq K_k\leq \infty$, then
$$\lambda+\mu+k\nu\geq -\frac{C_k}{t+\frac{1}{K_k}}.
$$
  \end{corollary} \begin{pf} For fixed $x_0\in M,$ since
the solution is smooth, there is a small $r_0>0$ such that whenever
$t\in [0,T],$ $x\in B_{t}(x_0,r_0)$, we have
$$
|Rm|(x,t)\leq r_0^{-2}.
$$
For the proof of (i), let $A\rightarrow \infty,$
$\delta\rightarrow0$ in the Proposition \ref{P2.1}, we get the
desired estimate. Case (ii) follows from Proposition \ref{P2.2} by
letting $A\rightarrow \infty$.
 \end{pf}
 \vskip 0.2cm
 In particular, in dimension 3, if the sectional curvature is nonnegative
 at $t=0,$ then this property is preserved for $t>0$ for \textbf{any} complete
 solutions.

Furthermore, for complete ancient solution,  for any fixed
$t\in(-\infty,0]$, by  Corollary \ref{C2.3} (ii) , we have
$(\lambda+\mu+k\nu )(t)\geq -\frac{C_k}{t-(-T)}$ for any $T>0.$
Since $C_k$ depends only on $k$, we have $(\lambda+\mu+k\nu)(t)\geq
0 $ for any $k\in \mathbb{Z}_{+}.$ This implies $\nu \geq 0,$ i.e.
the sectional curvature is nonnegative.

\begin{corollary}
\label{C2.4} Any ancient smooth complete  solution to the Ricci
flow (not necessarily having bounded curvature) on three manifold
must have nonnegative sectional curvature.
\end{corollary}
\begin{corollary}
\label{C2.5} Any ancient smooth complete  solution to the Ricci
flow (not necessarily having bounded curvature) must have
nonnegative scalar curvature.
\end{corollary}
\section{A priori estimates}
\subsection{}
 We will prove the following preliminary interior estimate, which holds for any dimension.

\begin{theorem}
\label{T3.1}
 There is a constant $C=C(n)$ with the following property.
 Suppose we have a smooth solution to the Ricci flow
 $(g_{ij})_{t}=-2R_{ij},$ $0\leq t\leq T,$  on an $n$-manifold $M$  such that
 $B_{t}(x_0,r_0),$ $0\leq t\leq T,$  is compactly contained in $M$ and

 (i)  $|Rm|\leq r_0^{-2}$ on $B_0(x_0,r_0)$ at $t=0$;

 (ii)  $$|Rm|(x,t)\leq
 \frac{K}{t}$$
 where $K\geq 1$,   $d_{t}(x,t)=dist_{t}(x_0,x)<
 r_0,$ whenever $0\leq t\leq T.$

 Then  we have
 $$|Rm|(x,t)\leq e^{CK}(r_0-d_{t}(x_0,x))^{-2}$$
whenever $0\leq t\leq T,$ $d_{t}(x,t)=dist_{t}(x_0,x)<r_0.$
 \end{theorem}

\begin{pf} By scaling, we may assume $r_0=1$ .

Since the result holds trivially by assumption when $t\geq 1.$
Without loss of generality, we may assume $T\leq 1.$

We argue by contradiction. Suppose we have a sequence of
$\delta\rightarrow0$, and a sequence of solutions satisfying the
assumptions in Theorem \ref{T3.1}. But
$|Rm|(x_1,t_1)>e^{\frac{K}{\delta}}\varepsilon^{-2}$ holds for
some point $(x_1,t_1)$, $d_{t_1}(x_1,x_0)<1-\varepsilon$, $t_1\in
[0,T].$

For any fixed $B\geq 1,$  by a point-picking technique of Perelman
\cite{P1}, we can choose another point $(\bar{x},\bar{t})$,
$\bar{x}\in B_{\bar{t}}(x_0,1-\frac{\varepsilon}{2} )$,
$\bar{t}\in (0,t_1]$ such that $\bar{Q}=|Rm|(\bar{x},\bar{t})\geq
e^{\frac{K}{\delta}}\varepsilon^{-2}$ and
 \be
 \label{3.1} |Rm|(x,t)\leq 2 \bar{Q}
\ee

whenever $d_{t}({x_0},x)\leq d_{\bar{t}}(\bar{x},x_0)+ 10 BK
\bar{Q}^{-\frac{1}{2}}$, $0\leq t\leq\bar{t}.$

 At the end of the
proof, it turns out that we only need to choose
$B=2\frac{e^{C(n)K}-1}{K}.$

Actually $(\bar{x},\bar{t})$ can be constructed as the limit of a
finite sequence $(x_i,t_i)$ satisfying  $0\leq t_k\leq t_{k-1}$,
$d_{t_k}(x_0,{x_k})\leq d_{t_{k-1}}(x_0,{x_{k-1}})+10BK
|Rm|({x_{k-1}},t_{k-1})^{-\frac{1}{2}}$, $|Rm|(x_k,t_k)\geq
2|Rm|(x_{k-1},t_{k-1})$. Since
 $$
|Rm|(x_k,t_k)\geq 2^{k-1}|Rm|(x_{1},t_{1})\geq
2^{k-1}e^{\frac{K}{\delta}}\varepsilon^{-2},
 $$
$d_{t_k}(x_0,{x_k})\leq
d_{t_1}(x_0,{x_{1}})+10BK\sum_{i=1}^{\infty}(2^{i-1}|Rm|(x_{1},t_{1}))^{-\frac{1}{2}}\leq
1-\varepsilon+40B Ke^{-\frac{K}{2\delta}}\varepsilon\leq
1-\frac{\varepsilon}{2}$. Clearly, if we choose
$B=2\frac{e^{C(n)K}-1}{K},$ the last inequality is guaranteed  by
 $e^{(C(n)-\frac{1}{2\delta})K}\leq \frac{1}{160},$ which holds trivially since $K\geq 1$ and $\delta\rightarrow 0.$
 Since the solution is smooth, this sequence must be finite
and the last element is what we want.

From this construction, we know $d_{\bar{t}}(\bar{x},x_0)+ 10 BK
\bar{Q}^{-\frac{1}{2}}\leq 1-\frac{\varepsilon}{2}.$

We denote by $C(n)$ various universal big constants depending only
upon the dimension. In the following argument, it may vary line by
line.

Now let $\varphi$ be a fixed smooth nonnegative non-increasing
cut-off function such that
 $\varphi=1$ on $(-\infty,d_{\bar{t}}(\bar{x},x_0)+BK
\bar{Q}^{-\frac{1}{2}}]$, $\varphi=0$ on
 $[d_{\bar{t}}(\bar{x},x_0)+ 10 BK
\bar{Q}^{-\frac{1}{2}},\infty).$ Clearly, we have
\be
 \label{3.2} |\varphi'|\leq C
\frac{\bar{Q}^{\frac{1}{2}}}{BK},
|\varphi''|+\frac{|\varphi'|^{2}}{\varphi}\leq
C\frac{\bar{Q}}{(BK)^{2}}. \ee

 Consider the function $u=\varphi({d_{t}({x}_0,x)})|Rm|(x,t)^{2},$
 it is clear
 \begin{equation*}
 \begin{split}
 (\frac{\partial}{\partial t}-\triangle)u& \leq\varphi'|Rm|^{2}(\frac{\partial}{\partial
t}-\triangle)d_{t}({x}_0,x)-2\varphi |\nabla Rm|^{2}\\&
-\varphi''|Rm|^{2}+C(n)\varphi |Rm|^{3}-2\nabla \varphi\cdot\nabla
|Rm|^{2}.
 \end{split}
 \end{equation*}
Since by (\ref{3.1}), $(\frac{\partial}{\partial
t}-\triangle)d_{t}({x}_0,x)\geq -C(n)\bar{Q}^{\frac{1}{2}}$
whenever $\bar{Q}^{-\frac{1}{2}} <d_{t}({x}_0,x).$ Then by the
maximum principle, and (\ref{3.1})(\ref{3.2}), it is clear that at
the maximum point,
\begin{equation*}
 \begin{split}
 \frac{d^{+}}{d t}u_{max}& \leq
 \frac{C(n)}{BK}|Rm|^{2}\bar{Q}
+C\varphi |Rm|^{3}\\
& \leq \frac{C(n)}{BK}\bar{Q}^{3}+C(n)\bar{Q}u_{max}(t).
 \end{split}
 \end{equation*}
Integrating this inequality, noting $u_{max}(0)\leq 1$ by
assumption, we get
$$
e^{-C(n)\bar{Q}t}u_{max}(t)\mid_{t=0}^{t=\bar{t}}\leq -
\frac{\bar{Q}^{2}}{BK}e^{-C(n)\bar{Q}t} \mid_{t=0}^{t=\bar{t}},
$$
and
$$
u_{max}(\bar{t})\leq
e^{C(n)\bar{Q}\bar{t}}+\frac{1}{BK}(e^{C(n)\bar{Q}\bar{t}}-1)\bar{Q}^{2}.
$$
Since $u_{max}(\bar{t})\geq u(\bar{x},\bar{t})=\bar{Q}^{2},$ and
$\bar{Q}\bar{t}\leq K,$ we have
$$
(1-\frac{e^{C(n)K}-1}{BK})\bar{Q}^{2}\leq e^{C(n)K}.
$$
Therefore,  if we choose $B=\frac{2(e^{C(n)K}-1)}{K}$, then we
have
$$
\bar{Q}\leq e^{C(n)K}
$$
which is a contradiction with $\bar{Q}\geq
e^{\frac{K}{\delta}}\varepsilon^{-2}$ as $\delta\rightarrow 0$.

This completes the proof of the Theorem \ref{T3.1}.
\end{pf}

\vskip 0.2cm
\begin{corollary}
\label{C3.2}
  Suppose we have a smooth solution to the Ricci flow
 $(g_{ij})_{t}=-2R_{ij},$ $0\leq t\leq T$,  such that
 that at $t=0$ we have $|Rm|\leq r_0^{-2}$ on $B_{0}(x_0,r_0)$; and  $$|Rm|(x,t)\leq \frac{K}{t}
 $$ whenever $0<t\leq T,$ $d_{0}(x,t)=dist_{0}(x_0,x)<
 r_0.$  Here we assume
 $B_{0}(x_0,r_0)$ is compactly contained in the manifold $M$.
 Then there is a  constant $C$ depending only on the dimension,
$$|Rm|(x,t)\leq e^{CK}(r_0-d_0(x_0,x))^{-2}$$
for $(x,t)\in B_{0}(x_0,r_0)\times [0,T].$

 \end{corollary}
\begin{pf} By \cite{P1}, for any fixed $p\in B_0(x_0,r_0)$,  as
long as the minimal geodesic $\gamma_t$ at time $t\in [0,r_0^{2}]$
connecting $p$ and $x_0$  lies in $B_0(x_0,r_0),$ we have
$$
\frac{d}{dt}d_{t}(x_0,p)\geq -C(n)\sqrt{\frac{K}{t}}.
$$
For any fixed $p\in B_0(x_0,r_0),$ let $[0,T')$ be the largest
interval such that any minimal geodesic $\gamma_t$ at time $t\in
[0,T']$ connecting $x_0$ and $p$ lies in $B_0(x_0,r_0)$ entirely. By
integrating the above inequality, we get
$$
d_{0}(x_0,p)\leq d_{t}(x_0,p)+C(n)\sqrt{K}\sqrt{T'}.
$$
This implies $B_t(x_0,\frac{r_0}{4})\subset
B_0(x_0,\frac{r_0}{2}),$ for any $t\in [0,\frac{r_0^{2}}{C(n)K}].$
By applying Theorem \ref{3.1} with
$T=\frac{r_0^{2}}{C(n)K}<(\frac{r_0}{4})^{2}$, there is a constant
$C(n)$ depending only on the dimension, such
 that
 $|Rm|\leq e^{C(n)K}r_0^{-2}$
whenever $0<t< \frac{r_0^{2}}{C(n)K},$ $d_{t}(x,t)=dist_{t}(x_0,x)<
\frac{1}{8} r_0.$ On the other hand, for $d_{0}(x_0,x)< r_0$ and
$t\in [\frac{r_0^{2}}{C(n)K},r_0^{2}]$, by assumption, we always
have
$$
|Rm|(x,t)\leq \frac{K}{t}\leq e^{C(n)K}r_0^{-2}.
$$
This in particular implies $|Rm|(x_0,t)\leq e^{C(n)K}r_0^{-2},$
for any $t\in [0,T].$

For any $x\in B_0(x_0,r_0),$ apply the above estimate on ball
$B_0(x, r_0-d_0(x_0,x))$ again, we know $|Rm|(x,t)\leq
e^{C(n)K}(r_0-d_0(x_0,x))^{-2}$ for any $t\in
[0,(r_0-d_0(x_0,x))^{2}].$ For $t>(r_0-d_0(x_0,x))^{2},$ we have
$|Rm|(x,t)\leq \frac{K}{t}\leq \frac{K}{(r_0-d_0(x_0,x))^{2}}\leq
e^{C(n)K}(r_0-d_0(x_0,x))^{-2}.$

The proof is completed.
\end{pf}

\subsection{}
We say a solution to the Ricci flow is ancient if it exists at
least on a half infinite time interval $(-\infty, T)$ for some
finite number $T.$ Ancient solution appears naturally in the blow
up  argument of singularity analysis of Ricci flow. The following
lemma will be used frequently  in the a priori estimates of this
section.
\begin{lemma}
 \label{L3.3}
 Let $g_{ij}(x,t),$ $t\in (-\infty,T)$ be a complete smooth
 non-flat ancient solution to the Ricci flow on an $n-$dimensional manifold
 $M$, with bounded and nonnegative curvature operator. Then for
  any $t\in (-\infty,T)$, the asymptotic volume ratio satisfies
 $$\nu_{M}(t):=\lim_{r\rightarrow
 \infty}\frac{vol_{t}(B_t(x,r))}{r^{n}}=0.$$
 \end{lemma}
This lemma was proved by \cite{P1}.
\begin{theorem}
 \label{T3.4}
 For any $C>0,$ there exists  $K>0$ with the following
 properties.
 Suppose we have a three dimensional smooth complete solution to the Ricci flow
 $(g_{ij})_{t}=-2R_{ij},$ $0\leq t\leq T,$ on a manifold $M,$ and assume
 that at $t=0$ we have $|Rm|(x,0)\leq r_0^{-2}$ on $B_0(x_0,r_0)$, and $R(x,0)\geq -r_0^{-2}$ on $M.$
  If
$g_{ij}(x,t)\geq
 \frac{1}{C}g_{ij}(x,0)$ for $x\in B_0(x_0,r_0)$, $t\in
 [0,r_0^{2}],$ then we have
 $$|Rm|(x,t)\leq 2
 r_0^{-2}$$
  whenever $0\leq t\leq \min\{ \frac{1}{K}r_0^{2},T\},$ $dist_{t}(x_0,x)< \frac{1}{K}r_0.$
 \end{theorem}
 \begin{pf} By scaling, let $r_0=1$.  By assumption $g_{ij}(x,t)\geq
 \frac{1}{C}g_{ij}(x,0)$, we have
 \begin{equation} \label{e3.2.0}
B_t(x_0,\frac{1}{\sqrt{C}})\subseteq B_0(x_0,1).
\end{equation}
 Let $T_0$ be the largest time such that $|Rm|(x,t)\leq 2$
 whenever $x\in B_{t}(x_0,\frac{1}{2\sqrt{C}}),$ $t\in [0,T_0].$
 We may assume $T_0<\min\{1,T\}.$ Otherwise, there is nothing to show.
 Hence there is a $(x_1,t_1)$ such that $|Rm|(x_1,t_1)=
 2,$ $x_1\in B_{t_1}(x_0,\frac{1}{2\sqrt{C}})$ and $t_1\leq T_0.$

  In the following arguments, we use
$\bar{C}$ to denote various constants depending only on $C$. By
using Corollary \ref{C2.3}, we know

$$
R(x,t)\geq -\bar{C},
$$
on $M\times [0,T_0].$  By  evolution equation of the volume
element $\frac{d}{dt}\log det(g)=-R$, this gives $\frac{det
(g)(t)}{det (g)(0)}\leq \bar{C}.$ Combining with the assumption
$g(t)\geq \frac{1}{C}g(0),$ we have
\begin{equation} \label{e3.2.1}
\frac{1}{\bar{C}}g(0)\leq g(t)\leq \bar{C}g(0)
\end{equation}
on $B_0(x_0,1)\times [0,T_0].$

Since the curvature on $B_{0}(x_0,1)$ of the initial metric $g$ is
bounded by $1$,  the exponential map (for the initial metric) at
$x_0$ is a local diffeomorphism from $B(0,1)\subset T_PM$ to the
geodesic ball $B_{0}(x_0,1),$ and such that $(\sin
1)\delta_{ij}\leq exp^{*}g_{ij}(x,0)\leq  (\sinh1)\delta_{ij}$ on
$B(0,1).$
 By the above estimate (\ref{e3.2.1}),  we have
\begin{equation} \label{e3.2.2}
\frac{1}{\bar{C}}\delta_{ij}\leq exp^{*}g_{ij}(x,t)\leq
\bar{C}\delta_{ij}
\end{equation}
on $B(0,1)\times [0,T_0].$ Let
$\bar{g}(\cdot,t)=\exp^{*}g(\cdot,t),$ then $\bar{g}(\cdot,t)$ is
a solution to the Ricci flow on the Euclidean ball $B(0,1)$,
moreover it is $\kappa$-noncollapsed for some $\kappa=\kappa(C)$
for all scales less than $1$ by (\ref{e3.2.2}).

Now we claim that there is a constant $K_0>0$ depending only on
${C}$ such that
\begin{equation} \label{e3.2.3}
|Rm|(x,t)\leq K_0
\end{equation}
as $x\in B_t(0,\frac{3}{4\sqrt{C}}),$ $t\in [0,T_0].$

Actually, suppose (\ref{e3.2.3}) is not true, then there is a
$(x_2,t_2)$ such that $|Rm|(x_2,t_2)\geq K_1\rightarrow \infty$,
$x_2\in B_{t_2}(0,\frac{3}{4\sqrt{C}}), 0<t_2\leq T_0.$ Now we can
choose another point $(\bar{x},\bar{t})$ so that
$\bar{Q}=|Rm|(\bar{x},\bar{t})\geq {K_1}$, $
\frac{1}{2\sqrt{C}}\leq d_{\bar{t}}(\bar{x},0)\leq
\frac{7}{8\sqrt{C}},$ $0<\bar{t}\leq t_2,$ and
\begin{equation}
\label{e3.2.4} |Rm|(x,t)\leq 4\bar{Q}
\end{equation}
for all $d_t(0,x)\leq
d_{\bar{t}}(0,\bar{x})+K_1^{\frac{1}{4}}\bar{Q}^{-\frac{1}{2}},$
$0\leq t\leq \bar{t}.$

Since $K\rightarrow \infty,$ we know
\begin{equation}
\label{e3.2.5}
d_{\bar{t}}(0,\bar{x})+K_1^{\frac{1}{4}}\bar{Q}^{-\frac{1}{2}}
\leq \frac{15}{16\sqrt{C}}.
\end{equation}
 Moreover by \cite{P1} and (\ref{e3.2.4}), it
follows $$\frac{d}{dt}d_{t}(0,\bar{x})\geq -\bar{C}\sqrt{\bar{Q}}$$
whenever $d_t(0,\bar{x})\leq
d_{\bar{t}}(0,\bar{x})+\frac{1}{2}K_1^{\frac{1}{4}}\bar{Q}^{-\frac{1}{2}}.$
By integrating this inequality, it is not hard to see
$d_t(0,\bar{x}) \leq
d_{\bar{t}}(0,\bar{x})+\bar{C}K_1^{\frac{1}{8}}\bar{Q}^{-\frac{1}{2}}$
whenever $0\leq \bar{Q}(\bar{t}-t)\leq
\min\{K^{\frac{1}{8}},\frac{\bar{Q}\bar{t}}{2}\}.$ Hence, if
$d_t(\bar{x},x)\leq K_1^{\frac{1}{8}}\bar{Q}^{-\frac{1}{2}},$ $0\leq
\bar{Q}(\bar{t}-t)\leq
\min\{K_1^{\frac{1}{8}},\frac{\bar{Q}\bar{t}}{2}\},$ we have
$d_t(0,x) \leq
d_{\bar{t}}(0,\bar{x})+\bar{C}K_1^{\frac{1}{8}}\bar{Q}^{-\frac{1}{2}}.$
By (\ref{e3.2.4}), this gives
\begin{equation}
\label{e3.2.6} |Rm|(x,t)\leq 4\bar{Q}, \ \ \ d_t(0,x)\leq
\frac{15}{\sqrt{16C}},
\end{equation}
for $x\in B_t(\bar{x},K^{\frac{1}{8}}\bar{Q}^{-\frac{1}{2}}),$ and
$0\leq \bar{Q}(\bar{t}-t)\leq
\min\{K^{\frac{1}{8}},\frac{\bar{Q}\bar{t}}{2}\}.$

 Recall in this
region, we always have (\ref{e3.2.2}) because of (\ref{e3.2.6})
and (\ref{e3.2.0}).

 Next, we will show \be
 \label{e}\bar{Q}\bar{t}\rightarrow
\infty, \ee which guarantees that the limit, which will be extracted
from a subsequence of the reascaled solutions around
$(\bar{x},\bar{t})$, is ancient.

Let $\varphi$ be a fixed smooth nonnegative non-increasing cut-off
function such that
 $\varphi=1$ on $(-\infty,d_{\bar{t}}(0,\bar{x})]$, $\varphi=0$ on
 $[d_{\bar{t}}(0,\bar{x})+K_1^{\frac{1}{4}}\bar{Q}^{-\frac{1}{2}},\infty).$

 Consider $u=\varphi(d_{t}(0,x))|Rm|(x,t)^{2},$ by applying the
 maximum principle as before, we have
\begin{equation*}
 \begin{split}
 \frac{d^{+}}{d t}u_{max}& \leq \bar{C} K_1^{-\frac{1}{4}} \bar{Q}^{3}+\bar{C}\bar{Q}u_{max}(t).
 \end{split}
 \end{equation*}
which gives
$$
\bar{Q}^{2}\leq e^{\bar{C}\bar{Q}\bar{t}}+
\bar{Q}^{2}\bar{C}K^{-\frac{1}{4}}(e^{\bar{C}\bar{Q}\bar{t}}-1).
$$
This implies $\bar{Q}\bar{t}\rightarrow \infty$ because
$\bar{Q}\geq {K_1}\rightarrow \infty.$

So by rescaling the solution around the point $(\bar{x},\bar{t})$
with the factor $\bar{Q}$ and shifting the time $\bar{t}$ to $0,$
and using Hamilton's compactness theorem and taking convergent
subsequence,  we get a smooth limit.   Note the curvature norm at
the new origin is $1$.  This limit is a nontrivial smooth complete
ancient solution to the Ricci flow with bounded curvature ($\leq
4$).  By Corollary \ref{C2.4}, this limit has nonnegative
curvature. But (\ref{e3.2.2}) indicates the asymptotic volume
ratio  of the limit is strictly positive, which is a contradiction
with Lemma \ref{L3.3}. So we have proved the claim (\ref{e3.2.3}).

 Let $\varphi$ be a fixed smooth nonnegative non-increasing
cut-off function such that
 $\varphi=1$ on $(-\infty,\frac{1}{2\sqrt{C}}]$, $\varphi=0$ on
 $[\frac{3}{4\sqrt{C}}, \infty).$ Consider the function
 $$u(x,t)=\varphi(d_{t}(0,x))|Rm|^{2}(x,t),$$ and by (\ref{e3.2.3}) and maximum
 principle, we obtain
 $$
 \frac{d^{+}}{dt}u_{\max}\leq \bar{C}
 $$
whenever $0\leq t\leq T_0.$  Recall we have $|Rm|(x_1,t_1)=2$ for
some $x_1\in B_{t_1}(x_0,\frac{1}{2\sqrt{C}})$ and $t_1\leq T_0.$
This gives $2\leq u_{\max}(t_1)\leq 1+\bar{C}t_1.$ Hence $T_0\geq
\frac{1}{\bar{C}}.$ The proof is completed.
\end{pf}
\begin{corollary}
\label{C3.5}

 For any $C,K_0>0$,  there exists a constant $K$ satisfying the following property.
 Suppose we have a three dimensional smooth complete solution to the Ricci flow
 $(g_{ij})_{t}=-2R_{ij},$ $0\leq t\leq T,$ on a manifold $M,$ and assume
 that at $t=0$ we have $|Rm|(\cdot,0)\leq K_0$ on $M$.
  If
$g_{ij}(\cdot,t)\geq
 \frac{1}{C}g_{ij}(\cdot,0)$ on $M\times
 [0,T],$ then  we have
 $$|Rm|(\cdot,t)\leq K $$
  for all $0\leq t \leq T.$
 \end{corollary}
\begin{pf} First of all, by Theorem \ref{T3.4} we know there is a
constant $T_0$ depending only on $K_0$ and $C$ such that
\begin{equation} \label{e3.2.8} |Rm|(\cdot,t)\leq 2K_0
\end{equation} for $0\leq t\leq \min\{T_0,T\}$ for some $T_0.$
Without loss of generality, we assume $T_0< T.$ By Corollary
\ref{C2.3} and assumption, we have
\begin{equation}
\label{e3.2.9} \frac{1}{\bar{C}}g(\cdot,0)\leq g(\cdot,t)\leq
\bar{C}g(\cdot,0)
\end{equation}
on $M\times [0,T].$

To prove the result, we will  argue by contradiction. Suppose there
is a point $(x_1,t_1)$ such that $|Rm|(x_1,t_1)\geq K\rightarrow
\infty.$ We can choose another point $(\bar{x},\bar{t})$ such that
$\bar{Q}=|Rm|(\bar{x},\bar{t})\geq {K},$ $\bar{t}\leq t_1$ and
$|Rm|(x,t)\leq 4\bar{Q}$, for all $d_t(x,\bar{x})\leq
K^{\frac{1}{4}}\bar{Q}^{-\frac{1}{2}}.$

Otherwise, we obtain a sequence of points $(x_k,t_k),$ such that
$t_1\geq t_2\geq\cdots,$ $|Rm|(x_k,t_k)\geq 4^{k-1}|Rm|(x_1,t_1),$
and $d_{t_k}(x_k,x_1)\leq \bar{C} K^{\frac{1}{4}} \sum
(4^{k-1}|Rm|(x_1,t_1))^{-\frac{1}{2}} \leq \bar{C}.$ Since
$d_{t_k}(x_k,x_1)\geq \frac{1}{\bar{C}}d_0(x_k,x_1),$ and the
solution is smooth, this procedure has to stop after a finite number
of steps. Now we pull back the solution locally by using the
exponential map (of the initial metric) at $\bar{x}$ to the
Euclidean ball of some fixed radius as before, and notice
$K^{\frac{1}{4}}\bar{Q}^{-\frac{1}{2}}\leq
\bar{C}K^{-\frac{1}{4}}\ll 1$ and (\ref{e3.2.9}). Then we can
rescale the solutions by the factor $\bar{Q}$ around
$(\bar{x},\bar{t})$ and extract a convergent subsequence. By
(\ref{e3.2.8}), the limit is ancient. The curvature (of the limit)
is bounded (by $4$). So by Corollary \ref{C2.4}, the limit has
nonnegative sectional curvature. It is clear by (\ref{e3.2.9}) and
the construction, the limit has maximal volume growth. So this is a
contradiction with Lemma \ref{L3.3}. The proof is completed.
\end{pf}
 \subsection{}
\begin{theorem}
\label{T3.6} For any $v_0>0,$ there is $K>0$ depending only on $v_0$
with the following properties. Let $(M,g(x,0))$ be a compete smooth
$3$-dimensional Riemannian manifold with nonnegative sectional
curvature,  $x_0\in M$ be a fixed point satisfying $|Rm|\leq
r_0^{-2}$ on $B_{0}(x_0,r_0)$ and $vol_0(B_0(x_0,r_0)))\geq
v_0r_0^{3},$ for some $r_0>0$.

Let $g(x,t)$, $t\in [0,T]$ be a smooth complete solution to the
Ricci flow with $g(x,0)$ as initial metric. Then we have
$$|Rm|(x,t)\leq 2r_0^{-2}$$ for all $x\in
B_{t}(x_0,\frac{r_0}{2})$, $0\leq t\leq
\min\{T,\frac{1}{K}r_0^{2}\}.$
\end{theorem}
\begin{pf} First of all, by Corollary \ref{C2.3}, for any $k>0$, there
is $C_k>0$ depending only on $k$ such that if at $t=0$,
$\lambda+\mu+k\nu\geq -K_k$ for some $0\leq K_k \leq\infty,$ then
for $t>0$, we have
$$\lambda+\mu+k\nu\geq -\frac{C_k}{t+\frac{1}{K_k}}.
$$
In our case, $\nu\geq0$ at $t=0,$ so we can choose $K_k=0$ for all
$k>0.$  Therefore, $\lambda+\mu+k\nu\geq 0$  for $t>0$, for any
$k>0.$ This implies  $\nu\geq0,$ i.e. curvatures are still
nonnegative for $t>0.$

By scaling, we assume $r_0=2.$

We imitate the proof of Theorem \ref{T3.4}. For the fixed $x_0\in
M,$ let $T_0$ be the largest time such that $|Rm(x,t)|\leq
\frac{1}{2}$ for all $x\in B_t(x_0,1)$ and $t\in [0,T_0].$ Recall by
assumption $|Rm(x,0)|\leq \frac{1}{4}$ on $B_0(x_0,2).$ Without loss
of generality, we assume $T_0<T.$ Then there is $(x_1,t_1)$ such
that $t_1\leq T_0$, $x_1\in B_{t_1}(x_0,1),$
$|Rm|(x_1,t_1)=\frac{1}{2}$. Our purpose is to estimate $T_0$ from
below by a positive constant depending only on $v_0.$

 Now we claim for fixed
$r>1$ there is a $B>0$ depending on $\frac{v_0}{r^{3}},$ such that
\begin{equation}
\label{e3.3.1} |Rm(x,t)|\leq B+B t^{-1}
\end{equation}
whenever $x\in B_{t}(x_0,\frac{r}{4})$ and $t\in [0,T_0].$

We will argue by contradiction. Actually, suppose (\ref{e3.3.1})
does not hold, then there is a sequence of solutions such that
there is some $(x_1,t_1),$ $x_1\in B_{t_1}(x_0,\frac{r}{4})$ and
$t_1\in [0,T_0]$ satisfying $|Rm(x_1,t_1)|\geq B+B t_1^{-1}$ with
$B\rightarrow \infty.$ By a point-picking technique of Perelman
\cite{P1}( Claim 1 and Claim 2 in Theorem 10 in \cite{P1}),
 we can
choose another $(\bar{x},\bar{t}),$ with
$\bar{Q}=|Rm|(\bar{x},\bar{t})\geq \frac{B}{\bar{t}}$ such that
\begin{equation}
\label{e3.3.2}
 |Rm|(x,t)\leq 4\bar{Q}
\end{equation}
for all $d_t(x,\bar{x})\leq
A^{\frac{1}{2}}\bar{Q}^{-\frac{1}{2}},$ $\bar{t}-A\bar{Q}^{-1}\leq
t\leq \bar{t},$ where $A$ tends to infinity with $B.$

Note that we have $vol_{t}(B_t(x_0,r))\geq \frac{v_0}{Cr^{3}}
r^{3},$ for $t\in [0,T_0].$ Since the curvature is nonnegative for
the solution, by volume comparison theorem, the solution is
$\kappa=\kappa(\frac{v_0}{r^{3}})$ non-collapsed on $B_{t}(x_0,r),$
for $t\leq T_0.$ So we can rescale the solution around
$(\bar{x},\bar{t})$ and extract a subsequence, finally obtain a
nontrivial ancient smooth complete solution to the Ricci flow, which
has maximal volume growth and bounded nonnegative curvature. This is
a contradiction with Lemma \ref{L3.3}. Therefore the claim
(\ref{e3.3.1}) is proved.

Now by choosing $r=8$ and applying Theorem \ref{T3.1},  we have
$|Rm(x,t)|\leq Const$ on $B_{t}(x_0,\frac{3}{2}),$ $t\in [0,T_0].$
Here the constant depends  only on $v_0.$

Consider the evolution equation of $\varphi(d_{t}(x_0,x))|Rm|(x,t),$
where $\varphi$ be a smooth nonnegative decreasing function which is
1 in $(-\infty,1]$ and $0$ in $[\frac{3}{2},\infty).$ As in the
proof of  Theorem \ref{T3.4}, by applying  maximum principle to the
equation of $\varphi(d_{t}(x_0,x))|Rm|(x,t),$  we conclude with
$T_0\geq \min\{T,\frac{1}{C} \}$. This completes the proof.
\end{pf}

\begin{corollary}
\label{C3.7} Let $(M,g(x))$ be a complete noncompact
$3$-dimensional manifold with bounded nonnegative sectional
curvature $0\leq Rm\leq K_0$,  for some fixed constants $K_0$. Let
${g}(x,t)$  be a smooth complete solution to the Ricci flow on
$M\times[0,T]$ with $g(x)$ as initial data. Then we have
$$0\leq Rm(\cdot,t)\leq
\frac{1}{\frac{1}{K_0}-{4t}}$$ for all $ 0\leq t <
\min\{T,\frac{1}{4K_0}\}$.
\end{corollary}
\begin{pf} First of all,
since we are considering the curvature estimate, by pulling back
the solution $g_t$ to the universal cover of the manifold, it is
sufficient to assume the manifold is simply-connected. We claim
that for such simply connected manifold, there is a constant
$i_0>0$ (may depend on the initial curvature bound) such that the
initial metric has injectivity radius bounded from below by
$i_0>0.$

Actually, since the curvature is bounded, by \cite{Sh1}, we may
deform the initial metric by the Ricci flow in a short time
interval $[0,\delta]$ such that the solution
$\tilde{g}_\tau$($\tau\in [0,\delta]$) has bounded and nonnegative
curvature $0\leq Rm\leq 2K_0$. Here the construction of
$\tilde{g}$ is from \cite{Sh1}, there should be no ambiguity with
the given solution in our theorem.  We have two possibilities. If
there is $\tau>0$ such that the sectional curvature vanishes
somewhere for some directions, then by the strong maximum
principle of Hamilton, the manifold splits as $\mathbb{R}\times
\Sigma$ or $\mathbb{R}^3$ metrically for all $\tau\in (0,\delta],$
where $\Sigma$ is surface with bounded and positive sectional
curvature. Another case is for all $\tau\in (0,\delta],$ the
sectional curvatures of the solution are positive everywhere. The
following fact is standard:  the injectivity radius of the
manifold is bounded from below by $\frac{\pi}{\sqrt{C}},$ for
simply connected closed even-dimensional manifold with $0<sec\leq
C$ or complete noncompact Riemannian manifold with  $0<sec\leq C.$
So we know the injectivity radius is bounded from below by a
uniform positive constant for any $0<\tau<\delta.$ Then our
assertion follows from the fact that the volume of the unit ball
at $t=0$ is uniformally bounded (from below by a positive
constant) by Ricci flow equation.

Note that by Corollary \ref{C2.3}(ii) the nonnegativity of
sectional curvature is preserved for $t>0.$ Since we have lower
injectivity radius bound at time $t=0$ from above argument, then
by applying Theorem \ref{T3.6}, we know there is a constant $K>0$
depending only on $i_0$ and $K_0$ such that $|Rm|(\cdot,t)\leq
2K_0$ for $ t\in [0,\min\{T,\frac{1}{K}\}].$ On the other hand,
once the curvature is bounded, we can apply the maximum principle
(on complete manifold with bounded curvature), yielding
\begin{equation}
\label{e3.3.3}0\leq Rm(x,t)\leq \frac{1}{\frac{1}{K_0}-4t}.
\end{equation}
Moreover we know the volume of the unit ball is also bounded from
below as long as the curvature is bounded.   So we may apply
Theorem \ref{T3.6} and maximum principle estimate repeatedly. So
(\ref{e3.3.3}) holds for all $ 0\leq t <\min\{T,\frac{1}{4K_0}\}$.
\end{pf}

 Combining \cite{CZ} and Corollary \ref{C3.7}, we complete the proof of Theorem \ref{T1.1}.
  The following
theorem follows also as a corollary of \cite{CZ} and Corollary
\ref{C3.5}.
\begin{theorem}
\label{T3.8}

Let $(M,g(0))$ be a complete smooth $3$ dimensional Riemannian
manifold such that $|Rm|(\cdot,0)\leq K_0$ on $M$.  Suppose we
have two  smooth complete solutions $g_1(t)$ and $g_2(t)$ to the
Ricci flow
 $(g_{ij})_{t}=-2R_{ij},$ $0\leq t\leq T,$ on $M$ with $g(0)$ as initial metric, and there is $C>0$
 such that
$g_{i}(\cdot,t)\geq
 \frac{1}{C}g(\cdot,0)$ on $M\times[0,T]$ $(i=1,2),$ then  we have
 $g_1(t)= g_2(t) $
  for all $0\leq t \leq T.$
 \end{theorem}

In concluding this section, we discuss the two dimensional case.
In this case, we can obtain purely local a priori estimates.

\begin{proposition}
\label{P3.9}  Let $g(x,t)$, $t\in [0,T]$ be a smooth solution to
the Ricci flow with $g(x,0)$ as initial metric on a two
dimensional Riemannian manifold $M$, $x_0\in M$. We assume
$B_{t}(x_0,r_0)$ is compactly contained in $M$ for any $t\in
[0,T];$ and at $t=0,$ $|R|(x,0)\leq r_0^{-2}$ on $B_{0}(x_0,r_0)$
and $vol_0(B_0(x_0,r_0))\geq v_0r_0^{2}$ for some constants
$r_0,v_0>0.$  Then there is a constant $K$ depending only on $v_0$
such that
$$|R|(x,t)\leq 2r_0^{-2}$$ for all $x\in
B_{t}(x_0,\frac{r_0}{2})$, $0\leq t\leq
\min\{T,\frac{1}{K}r_0^{2}\}.$
\end{proposition}
\begin{pf} The argument is similar to Theorem \ref{T3.6}. After
choosing the largest time $T_0$ such that curvature norm reaches
$2r_0^{-2}$ on the balls of radius $\frac{r_0}{2},$ by using
Proposition \ref{P2.1}, we have curvature estimate $R(x,t)\geq
-Cr_0^{-2}$ on balls of radius $\frac{3}{4}r_0$. Note that the
dimension is two, scalar curvature is the only curvature we have, so
this lower curvature bound enables us to apply the Bishop-Gromov
volume comparison theorem. Therefore, in the rest, we can argue as
in the proof of Theorem \ref{T3.6} to derive a lower bound for
$T_0.$
\end{pf}
The following result  is also clearly a corollary of Proposition
\ref{P3.9}.

\begin{theorem}
\label{T3.10}

Let $(M,g(0))$ be a complete smooth $2$ dimensional Riemannian
manifold such that $|R|\leq K_0,$   and $vol_0(B_0(\cdot,1))\geq
v_0$ for some fixed positive constants $K_0,v_0.$ Suppose we have
two smooth complete solutions $g_1(t)$ and $g_2(t)$ to the Ricci
flow
 $(g_{ij})_{t}=-2R_{ij},$ $0\leq t\leq T,$ on $M$ with  initial metric $g(0)$, then we have
 $g_1(t)= g_2(t) ,$
  for  $0\leq t \leq \min\{T,\frac{1}{K_0}\}.$
 \end{theorem}
\section{Concluding remarks}
It is interesting to know if the pseudolocality theorem of the
Ricci flow holds in a general class of Riemannian manifolds, and
the strong uniqueness theorem  holds in general as the corollary.
In particular, we may ask the question  for Euclidean space
$\mathbb{R}^{n}:$

\textbf{Question} \emph{ Does the strong uniqueness of the Ricci
flow hold on the Euclidean space $\mathbb{R}^{n}$ for $n\geq 4?$}

In the present paper, we have proved the case  for $n=2$ and $3.$

 We give
remarks for the analogous results on mean curvature flow. We
should mention that for codimension one hypersurfaces in Euclidean
space, the same type estimate was firstly established by Ecker and
Huisken \cite{EH}. There are much study for higher codimensional
mean curvature flow, see  M.T. Wang \cite{W}. A pseudolocality
estimate and  general strong
 uniqueness theorem for mean curvature flow were obtained in \cite{CY}.  The above question is  an intrinsic
 version of the result in \cite{CY}.


\begin{thebibliography}{99}
\bibitem{CaoZ} Cao, H.D. and Zhu,X.P., {\sl A complete proof
of Poincare and geometrization conjectures--application of
Hamilton-Perelman theory of Ricci flow}, Asian J. math. {\bf 10}2
(2006), 165-492.
\bibitem{CTY} Chau, A. Tam, L.F. and Yu, C.J., {\sl Pseudolocality for the Ricci flow and
applications,} arXiv: math.DG/0701153.
\bibitem{CY}  Chen, B. L. and Yin, L.,{\sl Uniqueness and pseudolocality theorems of the mean curvature
flow,} Comm. Anal. Geom. {\bf 15} (2007), no.3, 435-490.

\bibitem{CZ}  Chen, B. L. and Zhu, X. P.,
{\sl Uniqueness of the Ricci flow on complete noncompact
manifolds}, J. Diff. Geom. {\bf 74} (2006), 119-154.
\bibitem{CZ05F} Chen, B. L. and Zhu, X. P.,
 {\sl Ricci flow with Surgery on four-manifolds
 with positive isotropic
curvature}, J. Diff. Geom. {\bf 74} (2006), 177-264.
\bibitem{EH} Ecker, K. and Huisken, G. {\sl Interior estimates for hypersurfaces moving by mean curvature,}
Invent. Math. \textbf{105}, 547-569(1991).

\bibitem{Ha1}  Hamilton, R. S., {\sl Three manifolds with positive
Ricci curvature }, J. Diff. Geom. {\bf 17} (1982), 255-306.
\bibitem{Ha2}  Hamilton, R. S., {\sl The formation of singularities in
the Ricci flow}, Surveys in Differential Geometry (Cambridge, MA,
1993), {\bf 2}, 7-136, International Press, Combridge, MA,1995.
\bibitem{Ha3} Hamilton, R. S., {\sl Non-singular solutions to
the Ricci flow on three manifolds}, Comm. Anal. Geom. {\bf 1}
(1999), 695-729.
\bibitem{KL} Kleiner, B. and Lott,J., {\sl Note on Perelman's paper,}
http://www.math.lsa.umich.edu/research/ricciflow/perelman.html.
\bibitem{MT} Morgan,J.W. and Tian, G. {\sl Ricci flow and the
Poincare conjecture.} arXiv: math.DG/0607607.

 \bibitem{P1} Perelman, G., {\sl The entropy formula for the Ricci flow and its geometric
applications}, arXiv:math.DG/0211159 v1 November 11, 2002.
Preprint.
\bibitem{P2} Perelman, G., {\sl Ricci flow with surgery on three
manifolds} arXiv:math.DG/0303109 v1 March 10, 2003. Preprint.
\bibitem{Sh1}  Shi, W. X., {\sl Deforming the metric on complete
Riemannian manifold}, J. Diff. Geom. {\bf 30} (1989), 223-301.
\bibitem{Si} Simon,  M., {\sl Ricci flow of almost nonnegatively curved three
manifolds}, arXiv: math. DG/0612095 v1  Dec. 4, 2006. Preprint.
\bibitem{W}   Wang, M. T., {\sl  The mean curvature flow smoothes Lipschitz
submanifolds,} Comm. Anal. Geom. {\bf 12} (2004), no. 3, 581--599.

\end{thebibliography}
\end{document}